\documentclass[11pt]{amsart}

\usepackage{amsmath,amssymb,amscd,amsthm}

\makeatletter
\@namedef{subjclassname@2020}{\textup{2020} Mathematics Subject Classification}
\makeatother

\numberwithin{equation}{section}

\newcommand{\R}{\mathbb{R}}
\newcommand{\di}{\mathrm{d}}

\newcommand{\ddiv}{\mathrm{div}}

\newcommand{\trho}{\tilde{\rho}}
\newcommand{\tu}{\mathbf{\tilde{u}}}
\newcommand{\ess}{{\mathrm{ess}}}
\newcommand{\res}{{\mathrm{res}}}
\newcommand{\bfu}{{\bf u}}
\newcommand{\bfU}{{\bf U}}
\newcommand{\pat}{{\partial_t}}
\newcommand{\bfv}{{\bf v}}

\newtheorem{veta}{Theorem}[section]

{
\theoremstyle{definition}

\newtheorem{pozn}[veta]{Remark}
}

\begin{document}

\title[Relative energy inequality and weak-strong uniqueness]{Relative energy inequality and weak-strong uniqueness for an isothermal non-Newtonian compressible fluid}

\author{Richard Andr\'a\v sik, V\'aclav M\'acha, Rostislav Vod\'ak}

\begin{abstract}
Our paper deals with three-dimensional nonsteady Navier-Stokes equations for non-Newtonian compressible fluids. It contains a~derivation of the relative energy inequality for the weak solutions to these equations. We show that the standard energy inequality implies the relative energy inequality. Consequently, the relative energy inequality allows us to achieve a weak-strong uniqueness result. In other words, we present that the weak solution of the Navier-Stokes system coincides with the strong solution emanated from the same initial conditions as long as the weak solution exists.
\end{abstract}



\maketitle


\section{Introduction} 
\label{sec:intro}

This paper deals with the analysis of solution to the system describing a~compressible non-Newtonian fluid --  a fluid whose viscosity is non-constant, namely it depends on the shear rate in our case. 

The study of such system dates back to Mamontov \cite{Mam1}, \cite{Mam2} who showed the existence of a weak solution for an exponential growth of viscosity and an isothermal pressure. The existence of the weak solution in more general setting was provided by Zhikov and Pastukhova \cite{ZhPa}, nevertheless, their proof is wrong as noted by Feireisl, Liao and M\'alek in \cite{FeLiMa}. The article \cite{FeLiMa} itself deals with the existence of a weak solution, authors, however, had to assume a viscosity with rather an artificial term. 

The viscosity which does not possess an exponential growth admits only measure-valued solutions. This has been recently showed by Abbatiello, Feireisl and Novotn\'y in \cite{AbFeNo} for an isentropic pressure, and by Basari\'c in \cite{basaric} for the isothermal pressure. Authors of the former article also gave a proof of the relative energy and weak-strong uniqueness.

The existence of a strong solution was, up to now, showed in \cite{KaMaNe} where authors used the $L^p-L^q$ regularity approach to show the local-in-time well-posedness of the system describing compressible non-Newtonian fluids.

The relative entropy inequality (and corresponding weak-strong uniqueness) for the exponentially growing viscosity and isothermal pressure has not been investigated up to now. This is the main aim of our paper -- in particular, we derive a connection between results from \cite{Mam1}, \cite{Mam2} (resp. also \cite{basaric}) and \cite{KaMaNe}. 
\vspace{\baselineskip}

The motion of the non-Newtonian compressible fluid in a  domain $\Omega \subset \R^3$ is described by its velocity $\mathbf{u}:\Omega\to \mathbb R^3$ and density $\rho:\Omega\to [0,\infty)$. The time evolution of $\mathbf{u}$ and $\rho$ is governed by the continuity and momentum equations
\begin{eqnarray}
\partial_t \rho + \ddiv \left( \rho \mathbf{u} \right) & = & 0 ,  \label{eq:1} \\
\partial_t\left( \rho \mathbf{u} \right) + \ddiv\left( \rho \mathbf{u} \otimes \mathbf{u} \right) + \nabla p & = & \ddiv\ \mathbb{S} + \rho \mathbf{f} \quad \textrm{in } \Omega\times (0,T) , \label{eq:2}
\end{eqnarray}
where $T > 0$, $p$ stands for the pressure, $\mathbb{S}$ is the viscous stress tensor and $\mathbf{f}$ represents the external forces \cite{Mam1}.

Let us suppose that the fluid is isothermal and non-Newtonian. Namely, we assume that (without the loss of generality)
$$
\mathbb{S} = P(|D\mathbf{u}|)D\mathbf{u} ,\quad p = \rho
$$
where conditions on $P:[0,\infty)\to [0,\infty)$ are discussed later.

Equations~(\ref{eq:1}) and~(\ref{eq:2}) are supplemented with the no-slip boundary conditions
$$
\mathbf{u}|_{\partial \Omega} = 0 ,
$$
and initial conditions
$$
\rho(\cdot,0) = \rho_0\, ,\ (\rho\mathbf{u})(\cdot,0) = (\rho\mathbf{u})_0 \quad \mathrm{in}\ \Omega,
$$
where $\rho_0$ is non-negative.

The weak formulation of the system \eqref{eq:1}, \eqref{eq:2} endowed with presented boundary and initial condition is
\begin{equation}
\int_\Omega \rho_0 \varphi(\cdot,0)\ \di x + \int_0^T\int_\Omega \rho \partial_t \varphi + \rho \mathbf{u} \cdot \nabla\varphi\ \di x \di t = 0, \label{eq:var1}
\end{equation}
and
\begin{multline}
\int_\Omega (\rho\bfu)_0\psi(\cdot,0)\ \di x + \int_0^T\int_\Omega \left( \rho\mathbf{u} \cdot \partial_t\psi + (\rho\mathbf{u} \otimes \mathbf{u}) : D\psi + \rho\ddiv\psi \right)\ \di x \di t \\
= \int_0^T\int_\Omega \left( P(|D\mathbf{u}|)D\mathbf{u} : D\psi - \rho \mathbf{f}\cdot\psi \right)\ \di x \di t, \label{eq:var2}
\end{multline}
for any $\varphi \in \mathcal{C}_c^{\infty}(\bar{\Omega}\times\langle 0,T ))$ and $\psi \in  \mathcal{C}_c^{\infty}(\bar{\Omega}\times\langle 0,T ); \R^3)$.

The weak solution is supposed to satisfy the standard energy inequality expressed as \cite{Mam2} 
\begin{eqnarray}
\int_\Omega \left( \frac{1}{2}\rho|\mathbf{u}|^2 + \rho\ln{\rho}\right)(\cdot,\tau)\ \di x + \int_0^\tau\int_\Omega P(|D\mathbf{u}|)|D\mathbf{u}|^2\ \di x \di t \nonumber \\
\leq \int_0^\tau\int_\Omega \rho \mathbf{f} \cdot \mathbf{u}\ \di x \di t + \int_\Omega\left( \frac{|(\rho_0\mathbf{u})_0|^2}{2\rho_0} + \rho_0\ln{\rho_0} \right)\di x \label{eq:sei}
\end{eqnarray}
for almost all $\tau \in (0,T)$.

Let us define Young functions $\Phi_\gamma(z) = (1+z)\ln^\gamma{(1+z)}$, $\gamma > 1$. Functions $\Psi_\gamma$, denote their convex conjugates. For a given Young function $\Phi$ and its convex conjugate $\Psi$, we employ a standard notation of Orlicz class $\widetilde{L}_\Phi(\Omega)$ and Orlicz space $L_\Phi(\Omega)$. While $v \in \widetilde{L}_\Phi(\Omega)$ if $\int_\Omega \Phi(|v(x)|)\di x < +\infty$, $v \in L_\Phi(\Omega)$ if $\sup_w \int_\Omega |v(x)w(x)|\di x < +\infty$, where the supremum is taken over all functions $w \in \widetilde{L}_\Psi(\Omega)$ such that $\int_\Omega \Psi(|w(x)|)\di x \leq 1$. For further details about Orlicz spaces we refer to \cite{KuJoFu}.

Next, we define a Young function $M(z) = \mathrm{e}^z - z - 1$ and we denote by $N$ its convex conjugate. Similarly as in~\cite{Rosta2008}, we assume that the function $P$ satisfies the following five conditions\footnote{Hereinafter, we use the letter $C$ as an arbitrary positive constant which may vary from line to line, nevertheless, it is always independent of $\varrho$, $\bfu$, $U$ and $V$.}
\begin{equation} \label{eq:P1}
\int_{\Omega}P(|U|)|U|^2\ \di x \geq C\int_{\Omega}M(|U|)\ \di x,
\end{equation}
\begin{equation} \label{eq:P2}
\int_{\Omega}\left(P(|U|)U - P(|V|)V \right) : (U - V)\ \di x \geq 0,
\end{equation}
\begin{equation} \label{eq:P3}
P(z)|z|^2\ \textrm{is a convex function for } z \geq 0,
\end{equation}
\begin{equation} \label{eq:P4}
\int_{\Omega}N(P(|U|)|U|)\ \di x \leq C\left( 1 + \int_{\Omega}M(|U|)\ \di x \right),
\end{equation}
\begin{equation} \label{eq:P5}
P( | U - \lambda V | )( U - \lambda V ) \stackrel{M}{\rightharpoonup} P(|U|)U,\ \textrm{for } \lambda \rightarrow 0,
\end{equation}
for any $U$, $V$ belonging to Orlicz class $[\widetilde{L}_M(\Omega)]^{3\times 3}$. 

Throughout this paper, we assume even stricter condition then \eqref{eq:P2}. Namely, we assume the existence of $q\in (1,\infty)$ such that
\begin{equation}\label{eq:P6}
\int_\Omega (P(|U|)U - P(|V|)V): (U-V)\ \di x\geq C\int_\Omega |U-V|^q\ \di x
\end{equation}
for all $U$ and $V$ as above.

\vspace{\baselineskip}
We also define the following function spaces in accordance with \cite{Mam1}:
\begin{eqnarray*}
& \mathrm{Y} = \{ \mathbf{v} | D\mathbf{v} \in L_M(\Omega\times(0,T))^{3\times 3},\ \mathbf{v}|_{\partial\Omega\times(0,T)} = 0\} ,\\
& \mathrm{X} = \{ \mathbf{v} | D\mathbf{v} \in L_M(\Omega)^{3\times 3},\ \mathbf{v}|_{\partial\Omega} = 0\},\ \|\mathbf{v}\|_\mathrm{X} = \|D\mathbf{v}\|_{L_M(\Omega)}, \\
\end{eqnarray*}

First, our paper contains a derivation of the relative energy inequality for the weak solution constructed by Mamontov in \cite{Mam1} and \cite{Mam2}. This derivation is described in Section \ref{sec:rei}. The relative energy inequality is the cornerstone for further qualitative analysis of solutions as it allows to deduce various results concerning long-time behavior, singular limits, dimension reduction and weak-strong uniqueness result among others. The latter is performed in Section \ref{sec:wsu} -- we present that all weak solutions constructed by Mamontov are equal to the (unique) strong solution. Here, we would like to remind that the strong solution for the system in question with the periodic boundary conditions was constructed recently in \cite{KaMaNe}. Finally, Section \ref{sec:P} focuses on assumptions imposed on the function $P$. In particular, we provide a more convenient condition replacing (\ref{eq:P6}).

\vspace{\baselineskip}

\vspace{\baselineskip}
\section{Relative energy inequality} 
\label{sec:rei}

Let us define
$$
H(\rho) = \rho\int_1^\rho \frac{1}{z}\ \di z = \rho\ln{\rho}.
$$
Let us consider smooth functions $r$ and $\mathbf{U}$ such that $r$ is strictly positive and $\mathbf{U}$ satisfies the no-slip boundary conditions. Following \cite{FeireislREI}, {\em relative entropy} $\mathcal{E}([\rho,\mathbf{u}]|[r,\mathbf{U}])$ is defined as follows
\begin{eqnarray}
\mathcal{E}([\rho,\mathbf{u}]|[r,\mathbf{U}]) &=& \int_\Omega \left( \frac{1}{2}\rho|\mathbf{u} - \mathbf{U}|^2 + H(\rho) - H(r) - H'(r)(\rho - r) \right)\di x \nonumber \\
&=& \int_\Omega \left( \frac{1}{2}\rho|\mathbf{u} - \mathbf{U}|^2 + \rho\ln{\frac{\rho}{r} - (\rho - r)} \right)\di x. \label{eq:entropy}
\end{eqnarray}
We would like to point out that the function 
$$
\rho \mapsto H(\rho) - H(r) - H'(r)(\rho-r),\ \rho,\ r\geq 0.
$$
is strictly convex with minimum $0$ attained at $\rho = r$. Therefore, for every $0< \underline r<r<\overline{r}<\infty$ there is a positive constant $c$ such that 
\begin{equation}\label{coercivita.1}
H(\rho)- H(r) - H'(r)(\rho - r) > c(\rho - r)^2
\end{equation}
whenever $\rho \in \left(\frac{\underline r}{2},2\overline{r}\right)$, and 
\begin{equation}\label{coercivita.2}
H(\rho) - H(r) - H'(r)(\rho - r)> c|\rho -r|
\end{equation}
whenever $\rho \in \mathbb R^+ \setminus \left(\frac{\underline r}{2},2\overline{r}\right)$. 

Similarly as in \cite{FeireislSuitable}, a {\em suitable weak solution} to \eqref{eq:1} and \eqref{eq:2} is a couple $[\rho,\mathbf{u}]$ satisfying (\ref{eq:var1})-(\ref{eq:var2}), boundary and initial conditions, and the following relative energy inequality for all $r$ and $\bfU$ belonging to the class specified in Theorem \ref{veta:rei}:
\begin{eqnarray}
\mathcal{E}([\rho,\mathbf{u}]|[r,\mathbf{U}])(\tau) + \int_0^\tau \int_\Omega \left( P(|D\mathbf{u}|)D\mathbf{u} - P(|D\mathbf{U}|)D\mathbf{U} \right) : \left( D\mathbf{u} - D\mathbf{U} \right)\di x \di t \nonumber && \\
\leq \mathcal{E}([\rho_0,\mathbf{u}_0]|[r(\cdot,0),\mathbf{U}(\cdot,0)]) + \int_0^\tau \mathcal{R}(\rho,\mathbf{u},r,\mathbf{U})\ \di t, && \label{eq:rei}
\end{eqnarray}
where the remainder $\mathcal{R}$ is defined as
\begin{eqnarray}
\mathcal{R}(\rho,\mathbf{u},r,\mathbf{U}) &=& \int_\Omega \rho ( \partial_t \mathbf{U} + \mathbf{u} \nabla \mathbf{U}) \cdot (\mathbf{U} - \mathbf{u})\di x \nonumber \\
&& + \int_{\Omega} P(|D\mathbf{U}|)D\mathbf{U} : (D\mathbf{U} - D\mathbf{u})\di x + \int_\Omega \rho \mathbf{f} \cdot (\mathbf{u} - \mathbf{U})\di x \nonumber \\
&& + \int_\Omega \left( (r - \rho)\frac{\partial_t r}{r} + \frac{\nabla r}{r} \cdot (r\mathbf{U} - \rho\mathbf{u}) \right)\di x \nonumber \\
&& + \int_\Omega \ddiv\mathbf{U}(r - \rho)\di x. \label{eq:remainder}
\end{eqnarray}

If we set $r = r_0 = |\Omega|^{-1}\int_\Omega \rho_0\ \di x$ and $\mathbf{U} = 0$, then following the definition of the relative entropy (\ref{eq:entropy}), we arrive at
\begin{eqnarray*}
\mathcal{E}([\rho,\mathbf{u}]|[r_0,0]) &=& \int_\Omega \left( \frac{1}{2}\rho|\mathbf{u}|^2 + \rho\ln{\rho} - \rho\ln{r_0} \right)\di x, \\
\mathcal{E}([\rho_0,\mathbf{u}_0]|[r_0,0]) &=& \int_\Omega \left( \frac{1}{2}\rho_0|\mathbf{u}_0|^2 + \rho_0\ln{\rho_0} - \rho_0\ln{r_0} \right)\di x, \\
\mathcal{R}(\rho,\mathbf{u},r_0,0) &=& \int_\Omega \rho \mathbf{f} \cdot \mathbf{u}\ \di x.
\end{eqnarray*}
Since $\int_\Omega (\rho - \rho_0)\ln{r_0}\ \di x = \ln{r_0}\int_\Omega (\rho - \rho_0)\ \di x = 0$,  (\ref{eq:rei}) reduces to the standard energy inequality (\ref{eq:sei}). Thus, the relative energy inequality (\ref{eq:rei}) implies the standard energy inequality (\ref{eq:sei}). In section~\ref{sec:proof-rei}, we show that this implication holds also the other way round.

\begin{veta}\label{veta:rei}
Let $\rho$ and $\bfu$ be the weak solution specified in \eqref{eq:var1}, \eqref{eq:var2} and \eqref{eq:sei}. For some $p>1$ and $\gamma > 1$, let $\bfU\in Y$ be such that $D\bfU\in \widetilde L_{M}(\Omega \times (0,T))^9$ and $\pat \bfU\in L^p(0,T;L_{\Psi_\gamma}(\Omega))^3$. Let us further suppose that $\mathbf{f}\in L^p(0,T; L_{\Psi_\gamma}(\Omega))^3$ and $r:\Omega\times\langle0,T\rangle\to (0,\infty)$ is such that $r\in L^\infty(0,T;L_{\Phi_\gamma}(\Omega))$, $\pat \ln r\in L^1(0,T;L_{\Psi_\gamma}(\Omega))$ and $\nabla \ln r \in L^p(0,T;L_{\Psi_\gamma}(\Omega))^3$. Then the pairs $[\rho,\bfu]$ and $[r,\bfU]$ satisfy the relative entropy inequality \eqref{eq:rei}.
\end{veta}

\subsection{Proof of Theorem \ref{veta:rei}} \label{sec:proof-rei} 

First, we take $\frac{1}{2}|\mathbf{U}|^2 \psi$, $\psi=\chi_{\langle0,\tau\rangle}$ as a test function in (\ref{eq:var1}) -- this is possible although $\psi$ is not sufficiently regular as it suffices to take piecewise linear approximations of $\psi$ ensuring the validity of the following equalities for almost all $\tau\in(0,T)$. We infer
\begin{multline}
\frac{1}{2}\int_\Omega\rho_0 \bfU^2(\cdot,0)\ \di x + \int_0^\tau\int_\Omega \left( \rho\mathbf{U} \cdot \partial_t\mathbf{U} + \rho\mathbf{u} \nabla\mathbf{U} \cdot \mathbf{U} \right)\di x \di t\\ = \frac{1}{2}\int_\Omega \rho(\cdot,\tau)\bfU^2(\cdot,\tau)\ \di x, \label{eq:rei-step1}
\end{multline}
using identity $\partial_j(U_i U_i) = 2U_i\partial_j U_i$. Second, we test (\ref{eq:var2}) by $\mathbf{U}\psi$ (with the same remark as above) to arrive at
\begin{multline}
\int_\Omega (\rho\bfu)_0\bfU(\cdot,0)\, \di x +
 \int_0^\tau\int_\Omega ( \rho\mathbf{u} \cdot \partial_t\mathbf{U} + (\rho\mathbf{u}\otimes\mathbf{u}) : D\mathbf{U} + \rho\ddiv\mathbf{U}\\  + \rho \mathbf{f}\cdot\mathbf{U} - P(|D\mathbf{u}|)D\mathbf{u} : D\mathbf{U} )\, \di x \di t = \int_\Omega \rho\bfu(\cdot,\tau)\bfU(\cdot,\tau)\, \di x\label{eq:rei-step2}
\end{multline}
In the third step, we use $\ln{r}$ as a test function in (\ref{eq:var1}) to obtain
\begin{multline}
\int_\Omega \rho_0\ln{r}(\cdot,0)\, \di x + \int_0^\tau\int_\Omega \left( \rho\frac{\partial_t r}{r} + \rho\mathbf{u} \cdot \frac{\nabla r}{r} \right)\di x \di t\\ = \int_\Omega \rho(\cdot,\tau)\ln{r}(\cdot,\tau)\, \di x. \label{eq:rei-step3}
\end{multline}

We multiply (\ref{eq:rei-step1}) by $-1$ and sum it up with (\ref{eq:rei-step2}), (\ref{eq:rei-step3}) and the standard energy inequality (\ref{eq:sei}) to deduce
\begin{eqnarray}
&& \int_\Omega \left( \frac{1}{2}\rho|\mathbf{u}-\mathbf{U}|^2 + \rho\ln{\rho} - \rho\ln{r} - \rho \right)(\cdot,\tau)\di x \nonumber \\
&& + \int_0^\tau\int_\Omega \left( P(|D\mathbf{u}|)D\mathbf{u} - P(|D\mathbf{U}|)D\mathbf{U} \right) : \left( D\mathbf{u} - D\mathbf{U} \right)\di x \di t \nonumber \\
&& \leq \int_\Omega \left( \frac{1}{2}\rho_0|\mathbf{u}_0-\mathbf{U}(\cdot,0)|^2 + \rho_0\ln{\rho_0} - \rho_0\ln{r(\cdot,0)} - \rho_0 \right)\di x \nonumber \\
&& + \int_0^\tau\int_\Omega \rho\left( \partial_t\mathbf{U} + \mathbf{u}\nabla\mathbf{U} \right) \cdot (\mathbf{U} - \mathbf{u})\ \di x \di t + \int_0^\tau\int_\Omega \rho\mathbf{f}\cdot(\mathbf{u}-\mathbf{U})\ \di x \di t \nonumber \\
&& + \int_0^\tau\int_\Omega P(|D\mathbf{U}|)D\mathbf{U} : (D\mathbf{U} - D\mathbf{u})\ \di x \di t \nonumber \\
&& - \int_0^\tau\int_\Omega \left( \rho\frac{\partial_t r}{r} + \rho\mathbf{u}\cdot\frac{\nabla r}{r} + \rho\ddiv\mathbf{U} \right)\di x \di t. \label{eq:rei-step4}
\end{eqnarray}
We add equality $$\int_\Omega r(\cdot,\tau)\ \di x - \int_\Omega r(\cdot,0)\ \di x = \int_0^\tau\int_\Omega \partial_t r \ \di x \di t$$ to get
\begin{eqnarray}
\mathcal{E}([\rho,\mathbf{u}]|[r,\mathbf{U}]) + \int_0^\tau\int_\Omega \left( P(|D\mathbf{u}|)D\mathbf{u} - P(|D\mathbf{U}|)D\mathbf{U} \right) : \left( D\mathbf{u} - D\mathbf{U} \right)\di x \di t && \nonumber \\
\leq \mathcal{E}([\rho_0,\mathbf{u}_0]|[r(\cdot,0),\mathbf{U}(\cdot,0)]) + \int_0^\tau\int_\Omega \rho\left( \partial_t\mathbf{U} + \mathbf{u}\nabla\mathbf{U} \right) \cdot (\mathbf{U} - \mathbf{u})\ \di x \di t && \nonumber \\
+ \int_0^\tau\int_\Omega \rho\mathbf{f}\cdot(\mathbf{u}-\mathbf{U})\ \di x \di t + \int_0^\tau\int_\Omega P(|D\mathbf{U}|)D\mathbf{U}| : (D\mathbf{U} - D\mathbf{u})\ \di x \di t && \nonumber \\
+ \int_0^\tau\int_\Omega \left( (r-\rho)\frac{\partial_t r}{r} - \rho\mathbf{u}\cdot\frac{\nabla r}{r} - \rho\ddiv\mathbf{U} \right)\di x \di t. && \label{eq:rei-step5}
\end{eqnarray}
Finally, the boundary condition for $\bfU$ yields
$$
\int_\Omega (\mathbf{U}\cdot\nabla r + r\ddiv\mathbf{U})\ \di x = \int_\Omega \ddiv(r\mathbf{U})\ \di x = \int_{\partial\Omega} (r\mathbf{U})\cdot\mathbf{n}\ \di S = 0,
$$
and thus we can add the term $\int_\Omega (\mathbf{U}\cdot\nabla r + r\ddiv\mathbf{U})\ \di x$ to the right-hand side of (\ref{eq:rei-step5}) which gives the demanded inequality.

\vspace{\baselineskip}
\section{Weak-Strong uniqueness} 
\label{sec:wsu}

The relative entropy inequality can be used to prove the weak-strong uniqueness principle. In other words, it enables us to show that weak and strong solutions of (\ref{eq:1})-(\ref{eq:2}) with the same boundary and initial conditions coincide as long as the strong solution exists.

\begin{veta}\label{veta:wsu}
Let $\tilde \rho$, $\tilde \bfu$ be a strong solution to \eqref{eq:1} and \eqref{eq:2} satisfying the Dirichlet boundary condition and ${\tilde\rho}(x,t)\geq C,\ C\in \mathbb{R}^+$. If $P$ satisfies \eqref{eq:P1}--\eqref{eq:P6}, then every weak solution $\rho$, $\bfu$ emanating from the same initial data $\tilde \rho(\cdot,0)$ and $\tilde \bfu(\cdot,0)$ is equal to $\tilde \rho$ and $\tilde \bfu$.
\end{veta}

\begin{pozn} We assume the strong solution satisfies \eqref{eq:1}, \eqref{eq:2} pointwisely and all terms in this formulation are well defined.  In particular, $\tilde\varrho \in C^1(\overline{\Omega} \times \langle0,T\rangle)$ and $\tilde \bfu \in C^1(\overline{\Omega} \times \langle0,T\rangle)^3$ with $\nabla_x \bfu \in C(\overline{\Omega} \times \langle0,T\rangle)^{3\times 3}$. It is worth to mention that the proof presented below works also for weak solutions with sufficient regularity.
\end{pozn}

\subsection{Proof of Theorem \ref{veta:wsu}}  

Let us consider $r = \trho$ and $\mathbf{U} = \tu$, where $[\trho,\tu]$ is the strong solution of (\ref{eq:1})-(\ref{eq:2}) and $\trho \geq C$, $C \in \R^+$. The idea of the proof is to show that all terms in (\ref{eq:remainder}) can be bounded by the means of the left-hand side of (\ref{eq:rei}) in order to use a Gronwall type argument. Recall that the assumptions on the strong solution imply that there is not any vacuum region, i.e., it holds that $0<\underline r<\tilde\rho<\overline r<\infty$ for appropriately chosen constants $\underline r$ and $\overline r$. 

First, we introduce a decomposition of a general function $G = G(\rho)$ into the essential and residual part, namely, 
$$
G = G_\ess + G_\res
$$
where 
$$
G_\ess := \left\{
\begin{array}{l}
G\quad \mbox{on } \rho \in \left(\frac{1}2 \underline r, 2\overline r\right)\\
0\quad \mbox{otherwise. }
\end{array}
\right.
$$
Due to the convexity of $H$ one can deduce the following coercivity properties (see also \eqref{coercivita.1} and \eqref{coercivita.2}):
\begin{equation}\label{E.coercive}
\mathcal E([\rho,\bfu]|[\trho,\tu])\geq C\int_\Omega\left(\rho|\bfu - \tu|^2 + |\rho - \trho|^2_\ess + 1_\res + \rho_\res\right)\ \di x.
\end{equation}

According to (\ref{eq:remainder}), it holds that
\begin{eqnarray}
\mathcal{R}(\rho,\mathbf{u},\trho,\tu) &=& \int_\Omega \rho ( \partial_t\tu + \mathbf{u} \nabla\tu - \mathbf{f}) \cdot (\tu - \mathbf{u})\ \di x \nonumber \\
&& + \int_{\Omega} P(|D\tu|)D\tu : (D\tu - D\mathbf{u})\ \di x \nonumber \\
&& + \int_\Omega \left( (\trho - \rho)\frac{\partial_t \trho}{\trho} + (\trho\tu - \rho\mathbf{u}) \cdot \frac{\nabla \trho}{\trho} \right)\di x \nonumber \\
&& + \int_\Omega \ddiv\tu(\trho - \rho)\ \di x. \label{eq:wsu-remainder}
\end{eqnarray}

Since $[\trho,\tu]$ is a strong solution of (\ref{eq:1})-(\ref{eq:2}), we can rearrange the momentum equation (\ref{eq:2}) as follows:
\begin{equation*}
\frac{1}{\trho} \left( \tu\partial_t\trho + \trho \partial_t \tu + \ddiv(\trho\tu)\tu + \trho\tu\nabla\tu \right) - \mathbf{f}  =  \frac{1}{\trho}\ddiv\left( P(|D\tu|)D\tu \right) - \frac{\nabla\trho}{\trho}
\end{equation*}
which by the means of the continuity equation (\ref{eq:1}) reduces into
$$
\partial_t \tu + \tu\nabla\tu - \mathbf{f} = \frac{1}{\trho}\ddiv\left( P(|D\tu|)D\tu \right) - \frac{\nabla\trho}{\trho}.
$$
Hence, the first term in (\ref{eq:wsu-remainder}) can be rewritten as
\begin{eqnarray}
&& \int_\Omega \rho ( \partial_t\tu + \tu \nabla\tu - \mathbf{f}) \cdot (\tu - \mathbf{u})\ \di x + \int_\Omega \rho (\mathbf{u} - \tu) \nabla\tu \cdot (\tu - \mathbf{u})\ \di x \nonumber \\
&& = \int_\Omega \frac{\rho}{\trho} \ddiv\left( P(|D\tu|)D\tu \right) \cdot (\tu - \mathbf{u})\ \di x - \int_\Omega \frac{\rho\nabla\trho}{\trho} \cdot (\tu - \mathbf{u})\ \di x \nonumber \\
&& + \int_\Omega \rho (\mathbf{u} - \tu) \nabla\tu \cdot (\tu - \mathbf{u})\ \di x, \label{eq:wsu-step1}
\end{eqnarray}
which leads to (using the identity $\int_\Omega D\mathbf{v} : D\mathbf{w} \di x = - \int_\Omega \ddiv( D\mathbf{v}) \cdot \mathbf{w}\ \di x$):
\begin{eqnarray}
\mathcal{R}(\rho,\mathbf{u},\trho,\tu) &=& \int_\Omega \rho (\mathbf{u} - \tu) \nabla\tu \cdot (\tu - \mathbf{u})\ \di x \nonumber \\
&& + \int_{\Omega} \frac{1}{\trho}(\rho - \trho)\ddiv\left( P(|D\tu|)D\tu \right) \cdot (\tu - \mathbf{u})\ \di x \nonumber \\
&& + \int_\Omega (\trho - \rho) \left( \frac{\partial_t \trho}{\trho} + \tu \cdot \frac{\nabla\trho}{\trho} \right) \di x \nonumber \\
&& + \int_\Omega \ddiv\tu(\trho - \rho)\ \di x. \label{eq:wsu-step2}
\end{eqnarray}
With respect to the continuity equation (\ref{eq:1}), the last two terms in (\ref{eq:wsu-step2}) cancel each other out, because
\begin{eqnarray*}
\frac{\partial_t \trho}{\trho} + \tu \cdot \frac{\nabla\trho}{\trho} + \ddiv\tu &=& \frac{1}{\trho}\left( \partial_t \trho + \tu \cdot \nabla\trho + \trho\ddiv\tu \right) \\
&=& \frac{1}{\trho}\left( \partial_t \trho + \ddiv(\trho\tu) \right) = 0.
\end{eqnarray*}
Thus
\begin{eqnarray}
\mathcal{R}(\rho,\mathbf{u},\trho,\tu) &=& \int_\Omega \rho (\mathbf{u} - \tu) \nabla\tu \cdot (\tu - \mathbf{u})\ \di x \nonumber \\
&& + \int_{\Omega} \frac{1}{\trho}(\rho - \trho)\ddiv\left( P(|D\tu|)D\tu \right) \cdot (\tu - \mathbf{u})\ \di x. \label{eq:wsu-step3}
\end{eqnarray}

Concerning the first term of the remainder (\ref{eq:wsu-step3}), we have (recall that $\nabla\tilde u$ is a bounded function)
\begin{eqnarray}
\left| \int_\Omega \rho (\mathbf{u} - \tu) \nabla\tu \cdot (\tu - \mathbf{u})\ \di x \right| &\leq& C 
\int_\Omega \frac{1}{2}\rho|\tu - \mathbf{u}|^2 \di x \nonumber \\
&\leq& C 
\mathcal{E}([\rho,\mathbf{u}]|[\trho,\tu]).
\end{eqnarray}
Regarding the second term in (\ref{eq:wsu-step3}) we decompose it into two parts and we use the regularity of the strong solution to get
\begin{multline*}
\left|\int_\Omega \frac 1{\tilde \rho}(\rho - \tilde \rho) \ddiv(P(|D\tilde \bfu|)D\tilde \bfu)\cdot (\tilde \bfu - \bfu)\ \di x\right|\\
\leq C\int_\Omega |\tilde \rho - \rho|_{\ess}|\tilde \bfu - \bfu|\ \di x + C\int_\Omega |\tilde \rho - \rho|_{\res}|\tilde \bfu - \bfu|\ \di x =:I_1 + I_2.
\end{multline*}
Next,
\begin{multline*}
I_1\leq C\int_\Omega \frac 1{\sqrt \rho} |\tilde \rho - \rho|_\ess \sqrt \rho |\tilde \bfu - \bfu|\ \di x \\
\leq C\left(\int_\Omega |\tilde \rho - \rho|_\ess^2\ \di x + \int_\Omega \rho|\tilde \bfu - \bfu|^2\ \di x\right)\leq c\mathcal E([\rho,\bfu],[\tilde \rho,\tilde \bfu])
\end{multline*}
according to \eqref{E.coercive}. We split $I_2$ once again as
$$
I_2 = C\int_{\Omega\cap \{\rho>2\overline r\}} |\tilde\rho - \rho||\tilde \bfu - \bfu| \ \di x + C\int_{\Omega \cap \{\rho<\frac 12 \underline r\}} |\tilde \rho - \rho||\tilde \bfu  - \bfu|\ \di x
$$
where the first integral may be treated in the same way as $I_1$ and second integral is estimated with the help of the Korn inequality as follows
\begin{multline}\label{eq:korn}
\int_{\Omega\cap \{\rho<\frac 12 \underline \rho\}} |\tilde \rho - \rho||\tilde \bfu  - \bfu | \ \di x\leq C\int_\Omega |\tilde \bfu - \bfu|\ \di x\\
\leq C\int_\Omega 1^{q'}_\res \ \di x + \delta \int_\Omega |\tilde \bfu - \bfu |^q\ \di x\leq C \mathcal E([\rho,\bfu], [\tilde \rho,\tilde \bfu]) + \delta\int_\Omega |D\tilde \bfu- D\bfu|^q
\end{multline}
where $q' = \frac q{q-1}$ and $\delta>0$ might be as small as needed -- in particular, we choose $\delta$ such that (using \eqref{eq:P6})
$$
\delta \int_\Omega |D\tilde \bfu - D\bfu|^q\ \di x\leq \frac 12 \int_\Omega (P(|D\tilde \bfu|)D\tilde\bfu - P(|D\bfu|)D\bfu)(D\tilde \bfu - D \bfu)\ \di x
$$
so the last term of \eqref{eq:korn} can be absorbed in the left hand side of \eqref{eq:rei}. Summing up,  \eqref{eq:rei} yields
$$
\mathcal E([\rho,\bfu]|[\tilde\rho,\tilde\bfu])(\tau) \leq C\int_0^\tau \mathcal E([\rho,\bfu]|[\tilde\rho,\tilde\bfu])(s)\ \di s
$$
and, consequently, 
$$
\mathcal E([\rho,\bfu]|[\tilde\rho,\tilde\bfu])\equiv 0.
$$

\section{Assumptions on $P$} \label{sec:P}
The most restrictive assumption on $P$ is \eqref{eq:P6} which is not needed in the proof of the existence given by Mamontov. This assumption is discussed throughout this section. First, we show that a natural choice of $P$ satisfies the condition \eqref{eq:P6}. Next, we formulate a replacement assumption which is more convenient than \eqref{eq:P6}.
\subsection{$P(z) = \frac 1{z^2}M(z)$}
Let us consider a function
$$
P(z) = 
\left\{
\begin{tabular}{ll}
$\frac{M(z)}{z^2}$, & for $z \neq 0$, \\
$0$,              & for $z = 0$.
\end{tabular}
\right.
$$
It satisfies all conditions~(\ref{eq:P1})-(\ref{eq:P5}) and it satisfies also \eqref{eq:P6} with $q=3$. Indeed, recall that 
$$
P(z)z = \frac z2 + \frac {z^2}6 + \frac{z^3}{24} + \ldots = \sum_{i=1}^\infty \frac {z^i}{(i+1)!}
$$
and therefore we may write $P(z)z = F'(z) + G'(z)$ where $F(z) = \frac{z^2}{4}$ and 
$$
G(z) = \int_0^z \frac{s^2}{6} + \frac{s^3}{24} + \frac{s^4}{120}+\ldots \ \di s.
$$
Further, we have
\begin{equation*}
P(|D\bfu|)D\bfu = P(|D\bfu|)|D\bfu| \frac{D\bfu}{|D\bfu|} = \frac{\partial}{\partial D}\left(F(D\bfu) + G(D\bfu)\right).
\end{equation*}
Hence,
\begin{multline}
\left( P(|D\mathbf{u}|)D\mathbf{u} - P(|D\mathbf{v}|)D\mathbf{v} \right) : \left( D\mathbf{u} - D\mathbf{v} \right)  \\
= \left( \frac{\partial}{\partial D}F(|D\mathbf{u}|) - \frac{\partial}{\partial D}F(|D\mathbf{v}|) \right) : \left( D\mathbf{u} - D\mathbf{v} \right)\\
 + \left( \frac{\partial}{\partial D}G(|D\mathbf{u}|) - \frac{\partial}{\partial D}G(|D\mathbf{v}|) \right) : \left( D\mathbf{u} - D\mathbf{v} \right). \label{eq:wsu-step4}
\end{multline}
Since $F$ complies with the $\Delta_2$-condition, we can apply Lemma 21 from \cite{Diening} and deduce a lower bound for the first term in (\ref{eq:wsu-step4}) as follows:
\begin{multline}
 \left( \frac{\partial}{\partial D}F(|D\mathbf{u}|) - \frac{\partial}{\partial D}F(|D\mathbf{v}|) \right) : \left( D\mathbf{u} - D\mathbf{v} \right) \\
 \geq C F''(|D\mathbf{u}| + |D\mathbf{v}|) |D\mathbf{u} - D\mathbf{v}|^2 \geq C |D\mathbf{u} - D\mathbf{v}|^2 
\label{eq:wsu-step5}
\end{multline}
The second term in (\ref{eq:wsu-step4}) is non-negative. This follows since $D\bfu:D\bfv\leq \frac12|D\bfu| + \frac12|D\bfv|$ and thus
\begin{multline*}
\left( \frac{\partial}{\partial D}G(|D\mathbf{u}|) - \frac{\partial}{\partial D}G(|D\mathbf{v}|) \right) : \left( D\mathbf{u} - D\mathbf{v} \right)\\ =\left(\frac{G'(|D\bfu|)}{|D\bfu|} D\bfu - \frac{G'(|D\bfv|)}{|D\bfv|} D\bfv\right):(D\bfu - D\bfv)\\
=\frac{G'(|D\bfu|)}{|D\bfu|} D\bfu: D\bfu - \frac{G'(|D\bfu|)}{|D\bfu|} D\bfu:D\bfv\\
- \frac{G'(|D\bfv|)}{|D\bfv|} D\bfv:D\bfu + \frac{G'(|D\bfv|)}{|D\bfv|} D\bfv:D\bfv\\
\geq \frac 12 \frac{G'(|D\bfu|)}{|D\bfu|} |D\bfu|^2 - \frac 12 \frac{G'(|D\bfu|)}{|D\bfu|} |D\bfv|^2\\
-\frac12\frac{G'(|D\bfv|)}{|D\bfv|} |D\bfu|^2 +  \frac12\frac{G'(|D\bfv|)}{|D\bfv|} |D\bfv|^2\\
=\frac 12 \left(\frac{G'(|D\bfu|)}{|D\bfu|} - \frac{G'(|D\bfv|)}{|D\bfv|}\right)\left(|D\bfu|^2 - |D\bfv|^2\right) \geq 0
\end{multline*}
where we used the fact that $G'(z)/z$ is non-negative and increasing for $z\geq 0$.

\subsection{Further note on \eqref{eq:P6}}
In general, we claim that \eqref{eq:P6} is fulfilled whenever $P$ is a non-decreasing function satisfying 
\begin{multline}\label{eq:P6.a}
P'(z)\geq0\ \mbox{and there exists } c>0\mbox{ such that }P(z)\geq cz^\alpha\\ \mbox{for all }z\geq 0 \ \mbox{and some }\alpha>0\, .
\end{multline}

In order to prove that \eqref{eq:P6.a} implies \eqref{eq:P6} we start with a generalization of \cite[Sect I.4., Lemma 4.4, p. 14]{DiB}. The following holds for every $U,V\in \mathbb R^{3\times 3}$:
\begin{multline*}
(P(|U|)U-P(|V|)V):(U-V)\\
=\left(\int_0^1\frac{d}{ds}\left[ 
P(|sU+(1-s)V|)(sU+(1-s)V)\right]\ ds\right) :(U-V)\\ 
=\int_0^1P(|sU+(1-s)V|)|U-V|^2\ ds\\
+ \int_0^1\frac{P'(|sU+(1-s)V)|}{|
sU+(1-s)V|} | (sU+(1-s)V):( U-V) |^2\ ds\\
\geq \int_0^1P(|sU+(1-s)V|)|U-V|^2\ ds = (*).
\end{multline*}

Now, let us assume for a while that $|U|\geq |U-V|$. Then
\begin{multline*}
(*)= \int_0^1P(|U+(1-s)(V-U)|)|U-V|^2\ ds\\
\geq  \int_0^1P\left(|\ |U|-(1-s)|V-U|\ |\right)|U-V|^2\ ds\\
\geq \int _0^1P(s|U-V|)|U-V|^2\ ds\\
= |U-V|\int _0^{|U-V|}P(z)\ dz \geq c |U-V|^{2+\alpha}.
\end{multline*}

On the other hand, if $|U|<|U-V|$, we deduce that
\begin{equation*}
(*)= |U-V|^2 \int_0^1\frac{P(|sU+(1-s)V|)|sU+(1-s)V|^2}{|sU+(1-s)V|^2}\ ds =(1)\, . 
\end{equation*}
Since 
$$
|sU+(1-s)V|^2=|U+(1-s)(V-U)|^2\leq (2-s)^2|V-U|^2,
$$
we derive
\begin{multline*}
(1)\geq |U-V|^2 \int_0^1\frac{P(|sU+(1-s)V|)|sU+(1-s)V|^2}{(2-s)^2|U-V|^2}\ ds\\
\geq\frac{1}{4}\int_0^1P(|sU+(1-s)V|)|sU+(1-s)V|^2\ ds\\
\geq c\int_0^1Q(|sU+(1-s)V|^2)|sU+(1-s)V|^2\ ds=(2) ,
\end{multline*}
where $Q(z) = z^{\alpha/2}$. Since $G:= z^{1+\alpha/2}$ is a convex and non-decreasing function, we can use Jensen's inequality to deduce 
\begin{multline*}
(2)\geq C\, G\left(\int_0^1|sU+(1-s)V|^2\ ds\right)\\
= C\, G\left(\frac{1}{3}
(|U|^2+|V|^2+ U:V )\right)=(3).
\end{multline*}
Further,
$$
|U|^2 < |U-V|^2\Rightarrow 2 U:V < |V|^2,
$$
leads to 
\begin{multline*}
|U|^2+|V|^2+U:V\\
=\frac{9}{10}|U|^2+\frac{9}{10}|V|^2+\frac{6}{5} 
 U:V+\frac{1}{10}|U|^2+\frac{1}{10}|V|^2-\frac{2}{10} 
 U:V\\
\geq \frac{9}{10}|V|^2-\frac{6}{10}|V|^2+\frac{1}{10}(|U|^2+|V|^2-
2 U:V) \geq \frac{1}{10}|U-V|^2 .
\end{multline*}
Consequently,
$$
(3)\geq C\,G\left( c|U-V|^2\right) = C |U-V|^{2+\alpha}.
$$

We combine the two previous estimate to deduce that \eqref{eq:P6.a} implies \eqref{eq:P6} with $q = 2+\alpha$.


\vspace{\baselineskip}
\section*{Acknowledgments} 
This research was supported by The Ministry of Education, Youth and Sports CZ.02.1.01/0.0/0.0/17$\_$049/0008408 Hydrodynamic Design of Pumps.

The work of V\'aclav M\'acha was supported by Praemium Academi\ae \ of \v S. Ne\v casov\'a and by grant GA\v CR GA22-01591S in the framework of RVO:67985840.

\bibliography{references}

\providecommand{\bysame}{\leavevmode\hbox to3em{\hrulefill}\thinspace}
\providecommand{\MR}{\relax\ifhmode\unskip\space\fi MR }
\providecommand{\MRhref}[2]{%
  \href{http://www.ams.org/mathscinet-getitem?mr=#1}{#2}
}
\providecommand{\href}[2]{#2}
\begin{thebibliography}{10}

\bibitem{AbFeNo}
A.~Abbatiello, E.~Feireisl, and A.~Novotn\'{y}, \emph{Generalized solutions to
  models of compressible viscous fluids}, Discrete Contin. Dyn. Syst.
  \textbf{41} (2021), no.~1, 1--28. \MR{4182312}

\bibitem{basaric}
D.~Basari\'{c}, \emph{Existence of dissipative (and weak) solutions for models
  of general compressible viscous fluids with linear pressure}, J. Math. Fluid
  Mech. \textbf{24} (2022), no.~2, Paper No. 56, 22. \MR{4416230}

\bibitem{DiB}
E.~DiBenedetto, \emph{Degenerate parabolic equations}, Universitext,
  Springer-Verlag, New York, 1993. \MR{1230384}

\bibitem{Diening}
L.~Diening and F.~Ettwein, \emph{Fractional estimates for non-differentiable
  elliptic systems with general growth}, Forum Mathematicum \textbf{20} (2008),
  no.~3, 523--556.

\bibitem{FeireislREI}
E.~Feireisl, B.~J. Jin, and A.~Novotný, \emph{Relative entropies, suitable
  weak solutions, and weak-strong uniqueness for the compressible navier-stokes
  system}, J. Math. Fluid Mech. \textbf{14} (2012), 717--730.

\bibitem{FeLiMa}
E.~Feireisl, X.~Liao, and J.~M\'{a}lek, \emph{Global weak solutions to a class
  of non-{N}ewtonian compressible fluids}, Math. Methods Appl. Sci. \textbf{38}
  (2015), no.~16, 3482--3494. \MR{3423710}

\bibitem{FeireislSuitable}
E.~Feireisl, A.~Novotný, and Y.~Sun, \emph{Suitable weak solutions to the
  navier-stokes equations of compressible viscous fluids}, Indiana Univ. Math.
  J. \textbf{60} (2011), no.~2, 611--631.

\bibitem{KaMaNe}
M.~Kalousek, V.~M\'acha, and \v{S}. Ne\v{c}asov\'a, \emph{Local-in-time
  existence of strong solutions to a class of the compressible non-{N}ewtonian
  {N}avier-{S}tokes equations}, Mathematische Annalen (2021).

\bibitem{KuJoFu}
A.~Kufner, O.~John, and S.~Fu\v{c}\'{i}k, \emph{Function spaces}, Mechanics:
  Analysis, Springer Netherlands, 1977.

\bibitem{Mam1}
A.~E. Mamontov, \emph{Global solvability of the multidimensional navier-stokes
  equations of a compressible fluid with nonlinear viscosity i}, Siberian
  Mathematical Journal \textbf{40} (1999), no.~2, 351--362.

\bibitem{Mam2}
\bysame, \emph{Global solvability of the multidimensional navier-stokes
  equations of a compressible fluid with nonlinearly viscous fluid ii},
  Siberian Mathematical Journal \textbf{40} (1999), no.~3, 541--555.

\bibitem{Rosta2008}
R.~Vod\'{a}k, \emph{Asymptotic analysis of three dimensional navier-stokes
  equations for compressible nonlinearly viscous fluids}, Dynamics of PDE
  \textbf{5} (2008), no.~4, 299--311.

\bibitem{ZhPa}
V.~V. Zhikov and S.~E. Pastukhova, \emph{On the solvability of the
  {N}avier-{S}tokes system for a compressible non-{N}ewtonian fluid}, Dokl.
  Akad. Nauk \textbf{427} (2009), no.~3, 303--307. \MR{2573620}

\end{thebibliography}
\bibliographystyle{amsplain} 

\vspace{\baselineskip}
Richard Andr\'{a}\v{s}ik

\noindent
Department of Mathematical Analysis and Applications of Mathematics,
Faculty of Science, Palack\'y University Olomouc, 17. listopadu 12, 771 46 Olomouc, Czech Republic

\textit{E-mail address: andrasik.richard@gmail.com}

\vspace{\baselineskip}

V\'{a}clav M\'{a}cha

\noindent
Institute of Mathematics of the Czech Academy of Sciences, \v{Z}itn\'{a} 25, 115 67 Praha 1, Czech Republic

\textit{E-mail address: macha@math.cas.cz}

\vspace{\baselineskip}

Rostislav Vod\'{a}k

\noindent
Department of Mathematical Analysis and Applications of Mathematics,
Faculty of Science, Palack\'y University Olomouc, 17. listopadu 12, 771 46 Olomouc, Czech Republic

\textit{E-mail address: rostislav.vodak@gmail.com}

\end{document}